\documentclass{gtmon_a}
\pdfoutput=1

\usepackage{amscd}

\proceedingstitle{Proceedings of the Nishida Fest (Kinosaki 2003)}
\conferencestart{28 July 2003}
\conferenceend{8 August 2003}
\conferencename{International Conference in Homotopy Theory}
\conferencelocation{Kinosaki, Japan}

\editor{Matthew Ando}
\givenname{Matthew}
\surname{Ando}

\editor{Norihiko Minami}
\givenname{Norihiko}
\surname{Minami}

\editor{Jack Morava}
\givenname{Jack}
\surname{Morava}

\editor{W Stephen Wilson}
\givenname{W Stephen}
\surname{Wilson}

\title{Higher homotopy commutativity and cohomology of finite $H$--spaces}

\author{Yutaka Hemmi}
\givenname{Yutaka}
\surname{Hemmi}
\address{Department of Mathematics\\
Faculty of Science\\
Kochi University\\\newline
Kochi 780-8520\\
Japan}
\email{hemmi@math.kochi-u.ac.jp}
\urladdr{}

\author{Yusuke Kawamoto}
\givenname{Yusuke}
\surname{Kawamoto}
\address{Department of Mathematics\\
National Defence Academy\\\newline
Yokosuka 239-8686\\
Japan}
\email{yusuke@nda.ac.jp}
\urladdr{}

\volumenumber{10}
\issuenumber{}
\publicationyear{2007}
\papernumber{9}
\startpage{167}
\endpage{186}

\doi{}
\MR{}
\Zbl{}

\arxivreference{math.AT/0408215}

\keyword{higher homotopy commutativity}
\keyword{$H$--spaces}
\keyword{$A_p$--spaces}
\keyword{$AC_n$--spaces}
\keyword{Steenrod operations}
\subject{primary}{msc2000}{55P45}
\subject{primary}{msc2000}{57T25}
\subject{secondary}{msc2000}{55P48}
\subject{secondary}{msc2000}{55S05}

\received{14 August 2004}
\revised{19 September 2005}
\accepted{}
\published{29 January 2007}
\publishedonline{29 January 2007}
\proposed{}
\seconded{}
\corresponding{}
\version{}

\makeatletter
\def\cnewtheorem#1[#2]#3{\newtheorem{#1}{#3}[section]
\expandafter\let\csname c@#1\endcsname\c@theorem}

\AtBeginDocument{\let\bar\wbar\let\tilde\wtilde\let\hat\what}

\newtheorem{theorem}{Theorem}[section]
\cnewtheorem{lemma}[theorem]{Lemma}
\cnewtheorem{prop}[theorem]{Proposition}
\newtheorem{alphth}{Theorem}

\makeautorefname{alphth}{Theorem}

\theoremstyle{definition}
\cnewtheorem{definition}[theorem]{Definition}

\theoremstyle{remark}
\cnewtheorem{remark}[theorem]{Remark}

\makeatother  %  move after \newtheorem block

\numberwithin{equation}{section}

\begin{document}

\begin{htmlabstract}
We study connected mod p finite A<sub>p</sub>&ndash;spaces admitting
AC<sub>n</sub>&ndash;space structures with n&lt;p for an odd prime p.
Our result shows that if n&gt;(p-1)/2, then the mod p Steenrod algebra acts
on the mod p cohomology of such a space in a systematic way.  Moreover, we
consider A<sub>p</sub>&ndash;spaces which are mod p homotopy equivalent to
product spaces of odd dimensional spheres.  Then we determine the largest
integer n for which such a space admits an AC<sub>n</sub>&ndash;space
structure compatible with the A<sub>p</sub>&ndash;space structure.
\end{htmlabstract}

\begin{abstract}
We study connected mod $p$ finite $A_p$--spaces admitting $AC_n$--space
structures with $n<p$ for an odd prime $p$.  Our result shows that
if $n>(p-1)/2$, then the mod $p$ Steenrod algebra acts on the mod $p$
cohomology of such a space in a systematic way.  Moreover, we consider
$A_p$--spaces which are mod $p$ homotopy equivalent to product spaces
of odd dimensional spheres.  Then we determine the largest integer $n$
for which such a space admits an $AC_n$--space structure compatible with
the $A_p$--space structure.
\end{abstract}

\maketitle

\section{Introduction}
\label{sec1}
In this paper,
we assume that $p$ is a fixed odd prime
and that all spaces are localized at $p$
in the sense of Bousfield--Kan \cite{bousfield-kan}.

In the paper \cite{hemmi-kawamoto},
we introduced the concept of $AC_n$--space
which is an $A_n$--space whose multiplication
satisfies the higher homotopy commutativity
of the $n$-th order.
Then we showed that a mod~$p$ finite
$AC_n$--space with $n \ge p$ has the homotopy type of a torus.
Here by being mod $p$ finite,
we mean that the mod $p$ cohomology of the space is finite dimensional.
To prove it,
we first studied the action of the Steenrod operations
on the mod $p$ cohomology of such a space.
Then we showed that the possible cohomology generators
are concentrated in dimension $1$.

In the above argument,
the condition $n \ge p$ is essential.
In fact, any odd dimensional sphere admits an
$AC_{p-1}$--space structure
by \cite[Proposition~3.8]{hemmi-kawamoto}.
This implies that for any given exterior algebra,
we can construct a mod $p$ finite $AC_{p-1}$--space
such that the mod $p$ cohomology of it
is isomorphic to the algebra.

On the other hand,
if the $A_{p-1}$--space structure of the $AC_{p-1}$--space is
extendable to an $A_p$--space structure,
then the situation is different.
For example, it is known that an odd dimensional sphere
with an $A_p$--space structure
does not admit an $AC_{p-1}$--space structure
except for $S^1$.
In fact, an odd dimensional sphere $S^{2m-1}$
admits an $A_p$--space structure if and only if $m|(p-1)$,
and then it admits an $AC_n$--space structure
compatible with the $A_p$--space structure
if and only if $nm\le p$ by \cite[Theorem~2.4]{hemmi2}.
In particular, if $p=3$,
then mod~$3$ finite $A_3$--space with $AC_2$--space structure
means mod~$3$ finite homotopy associative
and homotopy commutative $H$--space.
Then by Lin \cite{lin5},
such a space has the homotopy type
of a product space of $S^1$s and $Sp(2)$s.

In this paper, we study mod~$p$ finite $A_p$--spaces with $AC_n$--space structures
for $n<p$.
First we consider the case of $n>(p-1)/2$.
In this case, we show the following fact
on the action of the Steenrod operations:

\begin{alphth}
\label{thmA}
Let $p$ be an odd prime.
If $X$ is a connected \textup{mod}~$p$ finite $A_p$--space
admitting an $AC_n$--space structure with $n>(p-1)/2$,
then we have the following:

\textup{(1)}\qua If $a\ge 0$,
$b>0$ and $0<c<p$,
then
\[
QH^{2p^a(pb+c)-1}(X;\mathbb{Z}/p)
=\mathscr{P}^{p^at}QH^{2p^a(p(b-t)+c+t)-1}(X;\mathbb{Z}/p)
\]
for $1\le t\le\min\{b,p-c\}$ and
\[
\mathscr{P}^{p^at}QH^{2p^a(pb+c)-1}(X;\mathbb{Z}/p)=0
\]
for $c\le t<p$.

\textup{(2)}\qua If $a\ge 0$ and $0<c<p$,
then
\[
\mathscr{P}^{p^at}\co QH^{2p^ac-1}(X;\mathbb{Z}/p)
\longrightarrow QH^{2p^a(tp+c-t)-1}(X;\mathbb{Z}/p)
\]
is an isomorphism for $1\le t<c$.
\end{alphth}

In the above theorem, the assumption $n>(p-1)/2$ is necessary.
In fact, (2) is not satisfied for the Lie group $S^3$ although
$S^3$ admits an $AC_{(p-1)/2}$--space structure
for any odd prime $p$ as is proved in \cite[Theorem~2.4]{hemmi2}.

\fullref{thmA} (1) has been already proved
for a special case or under additional hypotheses:
for $p=3$ by Hemmi \cite[Theorem~1.1]{hemmi1-1}
and for $p\ge5$ by Lin \cite[Theorem~B]{lin3}
under the hypotheses that the space
admits an $AC_{p-1}$--space structure
and the mod $p$ cohomology is $A_p$--primitively generated
(see Hemmi~\cite{hemmi3} and Lin~\cite{lin3}).

In the above theorem,
we assume that the prime $p$ is odd.
However, if we consider the case $p=2$,
then the condition $p>n>(p-1)/2$ is equivalent to $n=1$,
which means that the space is just an $H$--space.
Thus Theorem~A can be considered as the odd prime version
of Thomas \cite[Theorem~1.1]{thomas} or Lin \cite[Theorem~1]{lin1}.
(Note that in their theorems they assumed
that the mod $2$ cohomology of the space is primitively generated,
while we do not need such an assumption.)

By using \fullref{thmA},
we show the following result:

\begin{alphth}
\label{thmB}
Let $p$ be an odd prime.
If $X$ is a connected \textup{mod} $p$ finite $A_p$--space
admitting an $AC_n$--space structure
with $n>(p-1)/2$ and the Steenrod operations
$\mathscr{P}^j$ act on $QH^*(X;\mathbb{Z}/p)$ trivially
for $j\ge 1$,
then $X$ is \textup{mod} $p$ homotopy equivalent to a torus.
\end{alphth}

Next we consider the case of $n\le(p-1)/2$.
This includes the case $n=1$,
which means that the space is just
a mod $p$ finite $A_p$--space.
For the cohomology of mod $p$ finite $A_p$--spaces,
we can show similar facts to Theorem~A.
For example,
the results by Thomas \cite[Theorem~1.1]{thomas}
and Lin \cite[Theorem~1]{lin1} mentioned above
is for $p=2$,
and for odd prime $p$,
many results are known
(cf.\ \cite{akita},
\cite{hemmi1},
\cite{lin4}).

However,
for odd primes in particular,
those results have some ambiguities.
In fact,
there are many $A_p$--spaces with $AC_n$--space structures
for some $n\le(p-1)/2$ such that the Steenrod operations
act on the cohomology trivially.
In the next theorem,
we determine $n$ for which a product space
of odd dimensional spheres
to be an $A_p$--space with an $AC_n$--space structure.

\begin{alphth}
\label{thmC}
Let $X$ be a connected $A_p$--space
\textup{mod} $p$ homotopy equivalent
to a product space of odd dimensional spheres
$S^{2m_1-1}\times\cdots\times S^{2m_l-1}$
with $1\le m_1\le\cdots\le m_l$,
where $p$ is an odd prime.
Then $X$ admits an $AC_n$--space structure
if and only if $nm_l\le p$.
\end{alphth}

By the results of Clark--Ewing \cite{clark-ewing}
and Kumpel \cite{kumpel}, there are many spaces satisfying
the assumption of \fullref{thmC}.
Moreover, we note that the above result generalizes \cite[Theorem~2.4]{hemmi2}.

This paper is organized as follows:
In \fullref{sec2}, we first recall the modified projective
space $\mathscr{M}(X)$ of a finite $A_p$--space
constructed by Hemmi \cite{hemmi3}.
Based on the mod $p$ cohomology of $\mathscr{M}(X)$,
we construct an algebra $A^*(X)$
over the mod $p$ Steenrod algebra
which is a truncated polynomial algebra at height $p+1$
(\fullref{thm:modify}).
Next we introduce the concept of $\mathscr{D}_n$--algebra
and show that if $X$ is an $A_p$--space
with an $AC_n$--space structure,
then $A^*(X)$ is a $\mathscr{D}_n$--algebra (\fullref{thm:dn-algebra}).
Finally we prove the theorems in \fullref{sec3}
by studying the action of the Steenrod algebra
on $\mathscr{D}_n$--algebras algebraically
(\fullref{prop:propthmAB}
and \fullref{prop:propthmA}).

This paper is dedicated to Professor Goro Nishida
on his 60th birthday.
The authors appreciate the referee
for many useful comments.

\section{Modified projective spaces}
\label{sec2}
Stasheff \cite{stasheff1} introduced
the concept of $A_n$--space which is an $H$--space
with multiplication satisfying higher homotopy
associativity of the $n$-th order.
Let $X$ be a space and $n\ge 2$.
An $A_n$--form on $X$ is a family of maps
$\{M_i\co K_i\times X^i\to X\}_{2\le i\le n}$
with the conditions of \cite[I, Theorem~5]{stasheff1},
where $\{K_i\}_{i\ge 2}$ are polytopes called the associahedra.
A space $X$ having an $A_n$--form is called an $A_n$--space.
From the definition,
an $A_2$--space and an $A_3$--space
are the same as an $H$--space
and a homotopy associative $H$--space,
respectively.
Moreover,
it is known that an $A_{\infty}$--space
has the homotopy type of a loop space.

Let $X$ be an $A_n$--space.
Then by Stasheff \cite[I, Theorem~5]{stasheff1},
there is a family of spaces
$\{P_i(X)\}_{1\le i\le n}$
called the projective spaces
associated to the $A_n$--form on $X$.
From the construction of $P_i(X)$,
we have the inclusion
$\iota_{i-1}\co P_{i-1}(X)\to P_i(X)$ for $2\le i\le n$
and the projection $\rho_i\co P_i(X)\to P_i(X)/P_{i-1}(X) \simeq(\Sigma X)^{(i)}$
for $1\le i\le n$, where $Z^{(i)}$ denotes
the $i$-fold smash product
of a space $Z$ for $i\ge 1$.

For the rest of this section,
we assume that $X$ is a connected
$A_p$--space whose mod $p$ cohomology
$H^*(X;\mathbb{Z}/p)$ is an exterior algebra
\begin{equation}\label{eq:exterior}
H^*(X;\mathbb{Z}/p)\cong\Lambda(x_1,\ldots,x_l)
\qquad
\text{with $\deg x_i=2m_i-1$}
\end{equation}
for $1\le i\le l$,
where $1\le m_1\le\cdots\le m_l$.

Iwase \cite{iwase1} studied
the mod $p$ cohomology of the projective space
$P_n(X)$ for $1\le n\le p$.
If $1\le n\le p-1$, then there is an ideal $S_n\subset H^*(P_n(X);\mathbb{Z}/p)$
closed under the action of the mod $p$ Steenrod algebra
$\mathscr{A}_p^*$ such that
\begin{equation}\label{eq:cohomology of projective}
H^*(P_n(X);\mathbb{Z}/p)\cong T_n\oplus S_n
\qquad
\text{with $T_n=T^{[n+1]}[y_1,\ldots,y_l]$,
}
\end{equation}
where $T^{[n+1]}[y_1,\ldots,y_l]$ denotes
the truncated polynomial algebra at height $n+1$
generated by $y_i\in H^{2m_i}(P_n(X);\mathbb{Z}/p)$
with $\iota_1^*\ldots\iota_{n-1}^*(y_i)=\sigma(x_i)$
for $1\le i\le l$.
He also proved a similar result
for the mod $p$ cohomology of $P_p(X)$
under an additional assumption
that the generators $\{x_i\}_{1\le i\le l}$
are $A_p$--primitive
(see Hemmi~\cite{hemmi3} and Iwase~\cite{iwase1}).

Hemmi \cite{hemmi3} modified the construction
of the projective space $P_p(X)$ to get the algebra
$T^{[p+1]}[y_1,\ldots,y_l]$ also for $n=p$
without the assumption of the $A_p$--primitivity
of the generators.
Then he proved the following result:

\begin{theorem}[Hemmi {{\cite[Theorem~1.1]{hemmi3}}}]
\label{thm:modify}
Let $X$ be a simply connected $A_p$-space
whose \textup{mod} $p$ cohomology $H^*(X;\mathbb{Z}/p)$
is an exterior algebra in \eqref{eq:exterior},
where $p$ is an odd prime.
Then we have a space $\mathscr{M}(X)$
and a map $\epsilon\colon \Sigma X\to \mathscr{M}(X)$
with the following properties:

\textup{(1)}\qua
There is a subalgebra $R^*(X)$
of $H^*(\mathscr{M}(X);\mathbb{Z}/p)$ with
\begin{equation*}
R^*(X)\cong T^{[p+1]}[y_1,\ldots,y_l]\oplus M,
\end{equation*}
where $y_i\in H^{2m_i}(\mathscr{M}(X);\mathbb{Z}/p)$
are classes with $\epsilon^*(y_i)=\sigma(x_i)$
for $1\le i\le l$
and $M\subset H^*(\mathscr{M}(X);\mathbb{Z}/p)$
is an ideal with $\epsilon^*(M)=0$
and $M\cdot H^*(\mathscr{M}(X);\mathbb{Z}/p)=0$.

\textup{(2)}\qua
$R^*(X)$ and $M$ are closed under the action
of $\mathscr{A}_p^*$,
and so
\begin{equation}\label{eq:a(x)}
A^*(X)=R^*(X)/M\cong T^{[p+1]}[y_1,\ldots,y_l]
\end{equation}
is an unstable $\mathscr{A}_p^*$-algebra.

\textup{(3)}\qua
$\epsilon^*$ induces an $\mathscr{A}_p^*$-module
isomorphism:
\begin{equation*}
\begin{CD}
QA^*(X) @>>> QH^{*-1}(X;\mathbb{Z}/p).
\end{CD}
\end{equation*}
\end{theorem}

Next we recall the higher homotopy commutativity
of $H$--spaces.

Kapranov \cite{kapranov}
and Reiner--Ziegler \cite{reiner-ziegler}
constructed special polytopes $\{\Gamma_n\}_{n\ge 1}$
called the permuto--associahedra.
Let $n\ge 1$.
A partition of the sequence
$\mathbf{n}=(1,\ldots,n)$ of type $(t_1,\ldots,t_m)$
is an ordered sequence $(\alpha_1,\ldots,\alpha_m)$
consisting of disjoint subsequences $\alpha_i$
of $\mathbf{n}$ of length $t_i$
with $\alpha_1\cup\ldots\cup\alpha_m=\mathbf{n}$
(see Hemmi--Kawamoto~\cite{hemmi-kawamoto} and Ziegler~\cite{ziegler}
for the full details of the partitions).
By Ziegler~\cite[Definition~9.13, Example~9.14]{ziegler},
the permuto--associahedron $\Gamma_n$
is an $(n-1)$-dimensional polytope
whose facets (codimension one faces)
are represented by the partitions of $\mathbf{n}$
into at least two parts.
Let $\Gamma(\alpha_1,\ldots,\alpha_m)$
denote the facet of $\Gamma_n$
corresponding to a partition
$(\alpha_1,\ldots,\alpha_m)$.
Then the boundary
of $\Gamma_n$ is given by
\begin{equation}
\partial\Gamma_n=\bigcup_{(\alpha_1,\ldots,\alpha_m)}
\Gamma(\alpha_1,\ldots,\alpha_m)
\label{eqn:lambda}
\end{equation}
for all partitions $(\alpha_1,\ldots,\alpha_m)$
of $\mathbf{n}$ with $m\ge 2$.
If $(\alpha_1,\ldots,\alpha_m)$ is of type $(t_1,\ldots,t_m)$,
then the facet $\Gamma(\alpha_1,\ldots,\alpha_m)$
is homeomorphic to the product
$K_m\times\Gamma_{t_1}\times\cdots\times\Gamma_{t_m}$
by the face operator
$\epsilon^{(\alpha_1,\ldots,\alpha_m)}
\co K_m\times\Gamma_{t_1}\times\cdots\times
\Gamma_{t_m}\to\Gamma(\alpha_1,\ldots,\alpha_m)$
with the relations of \cite[Proposition~2.1]{hemmi-kawamoto}.
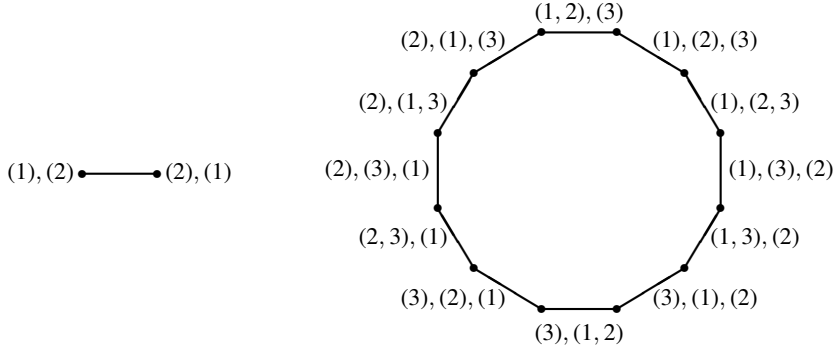
\begin{figure}[ht!]
\setlength{\unitlength}{1mm}
{\footnotesize
\begin{center}
\begin{picture}(20,50)(0,-23)

\put(0,0){\makebox(0,0)[r]{$(1),(2)$ }}
\put(10,0){\makebox(0,0)[l]{ $(2),(1)$}}

\thicklines
\put(0,0){\line(1,0){10}}
\put(0,0){\circle*{1.2}}
\put(10,0){\circle*{1.2}}
\end{picture}
\qquad\qquad\qquad\qquad
\begin{picture}(50,50)(-20,-5)
\put(0,38){\makebox(0,0)[b]{$(1,2),(3)$}}
\put(9.5,34.2){\makebox(0,0)[lb]{$(1),(2),(3)$}}
\put(16.4,27.4){\makebox(0,0)[l]{ $(1),(2,3)$}}
\put(18.8,18.4){\makebox(0,0)[l]{ $(1),(3),(2)$}}
\put(16.4,9.4){\makebox(0,0)[l]{ $(1,3),(2)$}}
\put(9.5,2.7){\makebox(0,0)[lt]{$(3),(1),(2)$}}
\put(0,-2){\makebox(0,0)[t]{$(3),(1,2)$}}

\put(-9.5,2.7){\makebox(0,0)[rt]{$(3),(2),(1)$}}
\put(-16.4,9.4){\makebox(0,0)[r]{$(2,3),(1)$ }}
\put(-18.8,18.4){\makebox(0,0)[r]{$(2),(3),(1)$ }}
\put(-16.4,27.4){\makebox(0,0)[r]{$(2),(1,3)$ }}
\put(-9.5,34.2){\makebox(0,0)[rb]{$(2),(1),(3)$}}

\thicklines
\put(-5,0){\line(1,0){10}}
\put(5,0){\line(5,3){9}}
\put(-5,0){\line(-5,3){9}}

\put(14,5.4){\line(3,5){4.8}}
\put(18.8,13.4){\line(0,1){10}}
\put(-14,5.4){\line(-3,5){4.8}}
\put(-18.8,13.4){\line(0,1){10}}

\put(18.8,23.4){\line(-3,5){4.8}}
\put(-18.8,23.4){\line(3,5){4.8}}

\put(-5,36.8){\line(1,0){10}}
\put(5,36.8){\line(5,-3){9}}
\put(-5,36.8){\line(-5,-3){9}}

\put(5,0){\circle*{1.2}}
\put(-5,0){\circle*{1.2}}
\put(14,5.4){\circle*{1.2}}
\put(-14,5.4){\circle*{1.2}}
\put(18.8,13.4){\circle*{1.2}}
\put(-18.8,13.4){\circle*{1.2}}
\put(18.8,23.4){\circle*{1.2}}
\put(-18.8,23.4){\circle*{1.2}}
\put(14,31.4){\circle*{1.2}}
\put(-14,31.4){\circle*{1.2}}
\put(5,36.8){\circle*{1.2}}
\put(-5,36.8){\circle*{1.2}}

\end{picture}
\end{center}
}
\caption{Permuto--associahedra $\Gamma_2$ and $\Gamma_3$}
\label{fig:gamma2gamma3}
\end{figure}
Moreover,
there is a family of degeneracy operators
$\{\delta_j\co\Gamma_i\to\Gamma_{i-1}\}_{1\le j\le i}$
with the conditions of \cite[Proposition~2.3]{hemmi-kawamoto}.

By using the permuto--associahedra,
Hemmi and Kawamoto \cite{hemmi-kawamoto}
introduced the concept of $AC_n$--form
on $A_n$--spaces.

Let $X$ be an $A_n$--space
whose $A_n$--form is given by
$\{M_i\}_{2\le i\le n}$.
An $AC_n$--form on $X$ is a family of maps
$\{Q_i\co\Gamma_i\times X^i\to X\}_{1\le i\le n}$
with the following conditions:
\begin{equation}
Q_1(*,x)=x.
\label{eqn:acnsp1}
\end{equation}
\begin{multline}
 Q_i(\epsilon^{(\alpha_1,\ldots,\alpha_m)}
(\sigma,\tau_1,\ldots,\tau_m),x_1,\ldots,x_i) \\
 = M_m(\sigma,Q_{t_1}(\tau_1,x_{\alpha_1(1)},\ldots,x_{\alpha_1(t_1)}),
\ldots,Q_{t_m}(\tau_m,x_{\alpha_m(1)},\ldots,x_{\alpha_m(t_m)}))
\label{eqn:acnsp2}
\end{multline}
for a partition $(\alpha_1,\ldots,\alpha_m)$
of $\mathbf{i}$ of type $(t_1,\ldots,t_m)$.
\begin{equation}
Q_i(\tau,x_1,\ldots,x_{j-1},*,x_{j+1},\ldots,x_i)
=Q_{i-1}(\delta_j(\tau),x_1,\ldots,x_{j-1},x_{j+1},\ldots,x_i)
\label{eqn:acnsp3}
\end{equation}
for $1\le j\le i$.

\begin{figure}[ht!]
\setlength{\unitlength}{1mm}
{\footnotesize
\begin{center}
\begin{picture}(20,50)(0,-23)

\put(0,0){\makebox(0,0)[r]{$xy$ }}
\put(10,0){\makebox(0,0)[l]{ $yx$}}

\thicklines
\put(0,0){\line(1,0){10}}
\put(0,0){\circle*{1.2}}
\put(10,0){\circle*{1.2}}
\end{picture}
\qquad\qquad\qquad\qquad
\begin{picture}(50,50)(-20,-5)
\put(2,37.8){\makebox(3,0)[lb]{$(xy)z$}}
\put(14,31.4){\makebox(0,0)[lb]{ $x(yz)$}}
\put(18.8,23.4){\makebox(0,0)[l]{ $x(zy)$}}
\put(18.8,13.4){\makebox(0,0)[l]{ $(xz)y$}}
\put(14,5.4){\makebox(0,0)[lt]{ $(zx)y$}}
\put(2,-1){\makebox(3,0)[lt]{$z(xy)$}}
\put(-5,-1){\makebox(3,0)[rt]{$z(yx)$}}
\put(-14.5,4){\makebox(0,0)[rt]{$(zy)x$ }}
\put(-18.8,13.4){\makebox(0,0)[r]{$(yz)x$ }}
\put(-18.8,23.4){\makebox(0,0)[r]{$y(zx)$ }}
\put(-14,31.4){\makebox(0,0)[rb]{$y(xz)$ }}
\put(-5,37.8){\makebox(3,0)[rb]{$(yx)z$}}

\thicklines
\put(-5,0){\line(1,0){10}}
\put(5,0){\line(5,3){9}}
\put(-5,0){\line(-5,3){9}}

\put(14,5.4){\line(3,5){4.8}}
\put(18.8,13.4){\line(0,1){10}}
\put(-14,5.4){\line(-3,5){4.8}}
\put(-18.8,13.4){\line(0,1){10}}

\put(18.8,23.4){\line(-3,5){4.8}}
\put(-18.8,23.4){\line(3,5){4.8}}

\put(-5,36.8){\line(1,0){10}}
\put(5,36.8){\line(5,-3){9}}
\put(-5,36.8){\line(-5,-3){9}}

\put(5,0){\circle*{1.2}}
\put(-5,0){\circle*{1.2}}
\put(14,5.4){\circle*{1.2}}
\put(-14,5.4){\circle*{1.2}}
\put(18.8,13.4){\circle*{1.2}}
\put(-18.8,13.4){\circle*{1.2}}
\put(18.8,23.4){\circle*{1.2}}
\put(-18.8,23.4){\circle*{1.2}}
\put(14,31.4){\circle*{1.2}}
\put(-14,31.4){\circle*{1.2}}
\put(5,36.8){\circle*{1.2}}
\put(-5,36.8){\circle*{1.2}}

\end{picture}
\end{center}
}
\caption{$Q_2(t,x,y)$ and $Q_3(\tau,x,y,z)$}
\label{fig:ac3form}
\end{figure}
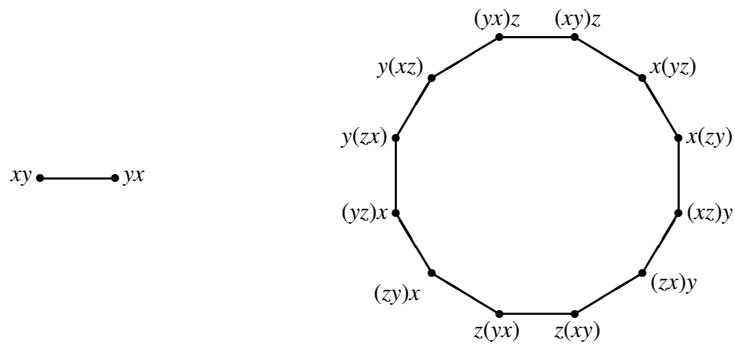

By \cite[Example~3.2 (1)]{hemmi-kawamoto},
an $AC_2$--form on an $A_2$--space
is the same as a homotopy commutative
$H$--space structure
since $Q_2\co\Gamma_2\times X^2\to X$
gives a commuting homotopy
between $xy$ and $yx$ for $x,y\in X$
(see \fullref{fig:ac3form}).
Let us explain an $AC_3$--form
on an $A_3$--space.
Assume that $X$ is an $A_3$--space
admitting an $AC_2$--form.
Then by using the associating homotopy
$M_3\co K_3\times X^3\to X$ and the commuting
homotopy $Q_2\co\Gamma_2\times X^2\to X$,
we can define a map
$\widetilde{Q}_3\co\partial\Gamma_3\times X^3\to X$
which is illustrated by the right dodecagon
in \fullref{fig:ac3form}.
For example, the uppermost edge represents the commuting
homotopy between $xy$ and $yx$,
and thus it is given by $Q_2(t,x,y)z$.
On the other hand, the next right edge is the associating
homotopy between $(xy)z$ and $x(yz)$
which is given by $M_3(t,x,y,z)$.
Then $X$ admits an $AC_3$--form
if and only if $\widetilde{Q}_3$ is extended
to a map $Q_3\co\Gamma_3\times X^3\to X$.
Moreover,
if $X$ is an $H$--space,
then by \cite[Example~3.2 (3)]{hemmi-kawamoto},
the multiplication of the loop space
$\Omega X$ on $X$ admits an $AC_{\infty}$--form.

Hemmi \cite{hemmi2} considered another concept of
higher homotopy commutativity of $H$--spaces.
Let $X$ be an $A_n$--space with the projective
spaces $\{P_i(X)\}_{1\le i\le n}$.
Let $J_i(\Sigma X)$ be the $i$-th stage
of the James reduced product space of $\Sigma X$
and $\pi_i\co J_i(\Sigma X)\to(\Sigma X)^{(i)}$
be the obvious projection for $1\le i\le n$.
A quasi $C_n$--form on $X$ is a family
of maps $\{\psi_i\co J_i(\Sigma X)\to P_i(X)\}_{1\le i\le n}$
with the following conditions:
\begin{align}
& \psi_1=1_{\Sigma X}\colon\Sigma X\longrightarrow\Sigma X.
&& \quad\quad\text{}
\label{eq:quasicn1} \\
& \psi_i|_{J_{i-1}(\Sigma X)}=\iota_{i-1}\psi_{i-1}
&& \quad\quad\text{for $2\le i\le n$.}
\label{eq:quasicn2} \\
& \rho_i\psi_i\simeq
\biggl(\sum_{\sigma\in\Sigma_i}\sigma\biggr)\pi_i
&& \quad\quad\text{for $1\le i\le n$,}
\label{eq:quasicn3}
\end{align}
where the symmetric group $\Sigma_i$
acts on $(\Sigma X)^{(i)}$
by the permutation of the coordinates
and the summation on the right hand side
is given by using the obvious co--$H$--structure
on $(\Sigma X)^{(i)}$ for $1\le i\le n$.

Hemmi and Kawamoto \cite{hemmi-kawamoto}
proved the following result:

\begin{theorem}[Hemmi--Kawamoto {\cite[Theorem~A]{hemmi-kawamoto}}]
\label{thm:acn-quasicn}
Let $X$ be an $A_n$--space for $n\ge 2$.
Then we have the following:

\textup{(1)}\qua
If $X$ admits an $AC_n$--form,
then $X$ admits a quasi $C_n$--form.

\textup{(2)}\qua
If $X$ is an $A_{n+1}$--space
admitting a quasi $C_n$--form,
then $X$ admits an $AC_n$--form.
\end{theorem}

\begin{remark}\label{rem:acn-quasicn}
In the proof of \fullref{thm:acn-quasicn} (2),
we do not need the condition \eqref{eq:quasicn3}.
In fact,
the proof of \fullref{thm:acn-quasicn} (2)
shows that if $X$ is an $A_{n+1}$-space
and there is a family of maps $\{\psi_i\}_{1\le i\le n}$
with the conditions \eqref{eq:quasicn1}--\eqref{eq:quasicn2},
then there is a family of maps $\{Q_i\}_{1\le i\le n}$
with the conditions \eqref{eqn:acnsp1}--\eqref{eqn:acnsp3}.
\end{remark}

Now we give the definition
of $\mathscr{D}_n$--algebra:

\begin{definition}\label{def:dn-algebra}
Assume that $A^*$ is an unstable $\mathscr{A}_p^*$--algebra
for a prime $p$.
Let $n\ge 1$.
Then $A^*$ is called a $\mathscr{D}_n$--algebra
if for any $\alpha_i\in A^*$ and $\theta_i\in\mathscr{A}_p^*$
for $1\le i\le q$ with
\begin{equation}\label{eq:dn-algebra-def1}
\sum_{i=1}^q\theta_i(\alpha_i)\in DA^*,
\end{equation}
there are decomposable classes
$\nu_i\in DA^*$ for $1\le i\le q$
with
\begin{equation}\label{eq:dnalgebra}
\sum_{i=1}^q\theta_i(\alpha_i-\nu_i)\in D^{n+1}A^*,
\end{equation}
where $DA^*$ and $D^tA^*$
denote the decomposable module and
the $t$-fold decomposable module of $A^*$
for $t>1$,
respectively.
\end{definition}

\begin{remark}\label{rem:dn-algebra}
It is clear from \fullref{def:dn-algebra}
that any unstable $\mathscr{A}_p^*$--algebra
is a $\mathscr{D}_1$--algebra.
On the other hand,
for an $A_p$--space $X$ which satisfies
the assumption of \fullref{thm:modify}
with $l\ge 1$,
the unstable $\mathscr{A}_p^*$--algebra
$A^*(X)$ given in \eqref{eq:a(x)}
cannot be a $\mathscr{D}_p$--algebra
since $\mathscr{P}^m(\alpha)=\alpha^p\ne 0$
for $\alpha\in QA^{2m}(X)$
from the unstable condition of $\mathscr{A}_p^*$
and $D^{p+1}A^*(X)=0$ in \eqref{eq:dnalgebra}.
\end{remark}

To prove \fullref{thmA} and \fullref{thmB},
we need the following theorem:

\begin{theorem}\label{thm:dn-algebra}
Let $p$ be an odd prime and $1\le n\le p-1$.
Assume that $X$ is a simply connected $A_p$-space
whose $\textup{mod}$ $p$ cohomology
$H^*(X;\mathbb{Z}/p)$ is an exterior algebra
in \eqref{eq:exterior}.
If the multiplication of $X$ admits a quasi $C_n$-form,
then $A^*(X)$ is a $\mathscr{D}_n$--algebra.
\end{theorem}

We need the following result
which is a generalization of Hemmi~\cite{hemmi2}:

\begin{lemma}\label{lem:da-algebra}
Assume that $X$ satisfies the same assumptions
as \fullref{thm:dn-algebra}.
If $\alpha_i\in H^*(P_n(X);\mathbb{Z}/p)$
and $\theta_i\in\mathscr{A}_p^*$ for $1\le i\le q$
satisfy
\begin{equation*}
\sum_{i=1}^q\theta_i(\alpha_i)=a+b
\qquad
\text{
with $a\in DH^*(P_n(X);\mathbb{Z}/p)$
and $b\in S_n$,}
\end{equation*}
then there are decomposable classes
$\nu_i\in DH^*(P_n(X);\mathbb{Z}/p)$ for $1\le i\le q$
with
\begin{equation*}
\sum_{i=1}^q\theta_i(\alpha_i-\nu_i)=b.
\end{equation*}
\end{lemma}

\begin{proof}
We give an outline of the proof
since the argument is similar to
Hemmi~\cite[Lemma~4.8]{hemmi2}.
It is clear for $n=1$.
If the result is proved for $n-1$,
then by the same reason as \cite[Lemma~4.8]{hemmi2},
we can assume $a\in D^nH^*(P_n(X);\mathbb{Z}/p)$.

Put
$\mathscr{U}_n=\widetilde{H}^*((\Sigma X)^{[n]};\mathbb{Z}/p)$,
$\mathscr{V}_n=QH^*(X;\mathbb{Z}/p)^{\otimes n}$
and
\begin{equation*}
\mathscr{W}_n=\bigoplus_{i=1}^n\widetilde{H}^*(X;\mathbb{Z}/p)^{\otimes i-1}
\otimes DH^*(X;\mathbb{Z}/p)\otimes
\widetilde{H}^*(X;\mathbb{Z}/p)^{\otimes n-i},
\end{equation*}
where $Z^{[n]}$ denotes the $n$-fold fat wedge
of a space $Z$ given by
\begin{equation*}
Z^{[n]}=\{(z_1,\ldots,z_n)\in Z^n\ |\
z_j=*\ \text{for some $1\le j\le n$}\}.
\end{equation*}
Then we have a splitting as an $\mathscr{A}_p^*$--module
\begin{equation}\label{split}
\widetilde{H}^*((\Sigma X)^n;\mathbb{Z}/p)
\cong\mathscr{U}_n\oplus\mathscr{V}_n
\oplus\mathscr{W}_n.
\end{equation}

Let $\mathscr{K}_n\co\widetilde{H}^*(X;\mathbb{Z}/p)^{\otimes n}
\to H^*(P_n(X);\mathbb{Z}/p)$
denote the following composite:
$$\widetilde{H}^*(X;\mathbb{Z}/p)^{\otimes n}
\xrightarrow{\sigma^{\otimes n}}
 H^*(\Sigma X;\mathbb{Z}/p)^{\otimes n} \\
\cong H^*((\Sigma X)^{(n)};\mathbb{Z}/p) \xrightarrow{\rho_n^*}
H^*(P_n(X);\mathbb{Z}/p).$$
Then by \cite[Theorem~3.5]{hemmi2},
there are $\widetilde{a}\in\mathscr{V}_n$
and $\widetilde{b}\in\mathscr{W}_n$
with $a=\mathscr{K}_n(\widetilde{a})$
and $b=\mathscr{K}_n(\widetilde{b})$.

Now we set $\lambda_n^*(\alpha_i)=c_i+d_i+e_i$
with respect to the splitting \eqref{split}
for $1\le i\le q$,
where $\lambda_n\co(\Sigma X)^n\to P_n(X)$
denotes the composite of $\psi_n$
with the obvious projection
$\omega_n\co(\Sigma X)^n\to J_n(\Sigma X)$.
From the same reason as Hemmi \cite[Lemma~4.8]{hemmi2}, we have
\begin{equation*}
\sum_{i=1}^q\theta_i(d_i)
=\sum_{\tau\in\Sigma_n}\tau(\widetilde{a})
=\lambda_n^*(a),
\end{equation*}
and so
\begin{equation*}
\lambda_n^*\biggl(\sum_{i=1}^q\theta_i(\mathscr{K}_n(d_i))\biggr)
=\sum_{\tau\in\Sigma_n}\tau\biggl(\sum_{i=1}^q\theta_i(d_i)\biggr)
=n!\sum_{\tau\in\Sigma_n}\tau(\widetilde{a})
=n!(\lambda_n^*(a)),
\end{equation*}
which implies
\begin{equation}\label{eq:a}
a=\frac{1}{n!}\sum_{i=1}^q\theta_i(\mathscr{K}_n(d_i))
\end{equation}
by \cite[Lemma~4.7]{hemmi2}.
If we put
\begin{equation*}
\nu_i=\frac{1}{n!}\mathscr{K}_n(d_i)
\in D^nH^*(P_n(X);\mathbb{Z}/p)
\end{equation*}
for $1\le i\le q$,
then by \eqref{eq:a},
$$\sum_{i=1}^q\theta_i(\alpha_i-\nu_i)=b,$$
which completes the proof.
\end{proof}

\begin{proof}[Proof of \fullref{thm:dn-algebra}]
From the construction of the space
$\mathscr{M}(X)$ in Hemmi~\cite[Section~2]{hemmi3},
we have a space $\mathscr{N}(X)$
and the following homotopy commutative diagram:
\begin{equation}\label{eq:dia1}
\begin{CD}
@. @. P_{p-2}(X) @>{\xi}>>
\mathscr{N}(X) @>{\eta}>> \mathscr{M}(X) \\
@. @. @| @V{\zeta}VV @VV{\kappa}V \\
P_n(X) @>>{\iota_n}> \cdots @>>{\iota_{p-3}}>
P_{p-2}(X) @>>{\iota_{p-2}}> P_{p-1}(X)
@>>{\iota_{p-1}}> P_p(X). \\
\end{CD}
\end{equation}
By \fullref{thm:modify} (2)
and \cite[page~593]{hemmi3},
we have that $M\subset R^*(X)$ is closed
under the action of $\mathscr{A}_p^*$
with $\eta^*(M)=0$,
which implies that $\eta^*|_{R^*(X)}\co R^*(X)
\to H^*(\mathscr{N}(X);\mathbb{Z}/p)$
induces an $\mathscr{A}_p^*$--homomorphism
$\mathscr{F}\co A^*(X)=R^*(X)/M\to
H^*(\mathscr{N}(X);\mathbb{Z}/p)$.
Then by applying the mod $p$ cohomology
to the diagram \eqref{eq:dia1},
we have the following commutative diagram
of unstable $\mathscr{A}_p^*$--algebras
and $\mathscr{A}_p^*$--homomorphisms:
\begin{equation*}
\begin{CD}
A^*(X) @>{\mathscr{F}}>> H^*(\mathscr{N}(X);\mathbb{Z}/p)
@<{\zeta^*}<< H^*(P_{p-1}(X);\mathbb{Z}/p) \\
@. @V{\xi^*}VV @VV{\iota_{p-2}^*}V \\
@. H^*(P_{p-2}(X);\mathbb{Z}/p)
@= H^*(P_{p-2}(X);\mathbb{Z}/p) \\
@. @. @VV{\iota_{p-3}^*}V \\
@. @. \vdots \\
@. @. @VV{\iota_n^*}V \\
@. @. H^*(P_n(X);\mathbb{Z}/p).
\end{CD}
\end{equation*}

First we assume $1\le n\le p-2$.
Put $\mathscr{G}_n(\alpha_i)=\beta_i$
for $1\le i\le q$,
where $\mathscr{G}_n\colon A^*(X)\to H^*(P_n(X);\mathbb{Z}/p)$
is the composite given by
$\mathscr{G}_n=\iota_n^*\ldots\iota_{p-3}^*\xi^*\mathscr{F}$.
Then by applying $\mathscr{G}_n$ to \eqref{eq:dn-algebra-def1},
we have
\begin{equation*}
\sum_{i=1}^q\theta_i(\beta_i)
\in DH^*(P_n(X);\mathbb{Z}/p),
\end{equation*}
and so by \fullref{lem:da-algebra},
there are decomposable classes
$\widetilde{\nu}_i\in DH^*(P_n(X);\mathbb{Z}/p)$
for $1\le i\le q$ with
\begin{equation}\label{eq:alpha-tilde-nu-tilde}
\sum_{i=1}^q\theta_i(\widetilde{\alpha}_i-\widetilde{\nu}_i)
=0.
\end{equation}
If we choose decomposable classes
$\nu_i\in DA^*(X)$
to satisfy $\mathscr{G}_n(\nu_i)=\widetilde{\nu}_i$
for $1\le i\le q$,
then by \eqref{eq:alpha-tilde-nu-tilde},
\begin{equation*}
\sum_{i=1}^q\theta_i(\alpha_i-\nu_i)
\in D^{n+1}A^*(X),
\end{equation*}
which completes the proof
in the case of $1\le n\le p-2$.

Next let us consider the case of $n=p-1$.
Put $\mathscr{F}(\alpha_i)=\widetilde{\alpha}_i
\in H^*(\mathscr{N}(X);\mathbb{Z}/p)$
for $1\le i\le q$.
Then we have
$$\sum_{i=1}^q\theta_i(\widetilde{\alpha}_i)
\in DH^*(\mathscr{N}(X);\mathbb{Z}/p).$$
By \cite[Proposition~5.2]{hemmi3},
we see that $\mathscr{F}(A^*(X))$ is contained
in $\zeta^*(H^*(P_{p-1}(X);\mathbb{Z}/p))$,
and so we can choose $\beta_i\in H^*(P_{p-1}(X);\mathbb{Z}/p)$
and $a\in DH^*(P_{p-1}(X);\mathbb{Z}/p)$
with $\zeta^*(\beta_i)=\widetilde{\alpha}_i$
and
\begin{equation*}
\zeta^*(a)=\sum_{i=1}^q\theta_i(\widetilde{\alpha}_i)
\end{equation*}
for $1\le i\le q$.
Then we can set
\begin{equation*}
\sum_{i=1}^q\theta_i(\beta_i)=a+b
\end{equation*}
with $\zeta^*(b)=0$,
and by \cite[Lemma~5.1]{hemmi3},
we have $b\in S_{p-1}$.
By \fullref{lem:da-algebra},
there are decomposable classes
$\mu_i\in DH^*(P_{p-1}(X);\mathbb{Z}/p)$
for $1\le i\le q$ with
\begin{equation*}
\sum_{i=1}^q\theta_i(\beta_i-\mu_i)=b.
\end{equation*}
Let $\nu_i\in DA^*(X)$
with $\mathscr{F}(\nu_i)=\zeta^*(\mu_i)$
for $1\le i\le q$.
Then we have
\begin{equation*}
\sum_{i=1}^q\theta_i(\alpha_i-\nu_i)\in D^pA^*(X),
\end{equation*}
which implies the required conclusion.
This completes the proof of \fullref{thm:dn-algebra}.
\end{proof}

\section[Proofs of \ref{thmA} and \ref{thmB}]{Proofs of \fullref{thmA}
and \fullref{thmB}}
\label{sec3}
In this section,
we assume that $A^*$ is an unstable
$\mathscr{A}_p^*$--algebra which is
the truncated polynomial algebra
at height $p+1$ given by
\begin{equation}\label{eq:truncated}
A^*=T^{[p+1]}[y_1,\ldots,y_l]
\qquad
\text{with $\deg y_i=2m_i$}
\end{equation}
for $1\le i\le l$,
where $1\le m_1\le\cdots\le m_l$.
Moreover, we choose the generators $\{y_i\}$ to satisfy
\begin{equation}\label{eq:generators}
\mathscr{P}^1(y_i)\in DA^* \text{ or }\mathscr{P}^1(y_i)=y_j
\qquad
\text{for some $1\le j\le l$.}
\end{equation}

The above is possible by the same argument
as Hemmi \cite[Section~4]{hemmi1}.

First we prove the following result:

\begin{prop}\label{prop:propthmAB}
Suppose that $A^*$ is a $\mathscr{D}_n$--algebra
and $1\le i\le l$.
If $\mathscr{P}^1(y_i)$ contains the term
$y_j^t$ for some $1\le j\le l$ and $1\le t\le n$,
then $y_j=\mathscr{P}^1(y_k)$ for some $1\le k\le l$.
\end{prop}

\begin{proof}
If $t=1$, then by \eqref{eq:generators},
the result is clear.
Let $t$ be the smallest integer
with $1<t\le n$ such that
the term $y_j^t$ is contained
in $\mathscr{P}^1(y_{i'})$ for some $1\le{i'}\le l$.
Then by \eqref{eq:generators},
we have $\mathscr{P}^1(y_{i'})\in DA^*$.
Since $A^*$ is a $\mathscr{D}_n$--algebra,
there is a decomposable class $\nu\in DA^*$
with $\mathscr{P}^1(y_{i'}-\nu)\in D^{n+1}A^*$.
This implies that $\mathscr{P}^1(\nu)$
contains the term $y_j^t$,
and so there is one of the generators
$y_{i''}$ of \eqref{eq:truncated} for $1\le{i''}\le l$
such that $\mathscr{P}^1(y_{i''})$
contains the term $y_j^s$
for some $1\le s<t$.
Then we have a contradiction,
and so $t=1$.
This completes the proof.
\end{proof}

In the proof of \fullref{thmA},
we need the following result:

\begin{prop}\label{prop:propthmA}
Let $p$ be an odd prime.
If $A^*$ is a $\mathscr{D}_n$--algebra with $n>(p-1)/2$,
then the indecomposable module $QA^*$
of $A^*$ satisfies the following:

\textup{(1)}\qua If $a\ge 0$,
$b>0$ and $0<c<p$,
then
\begin{equation}\label{eq:A1}
QA^{2p^a(pb+c)}={\mathscr{P}}^{p^at}QA^{2p^a(p(b-t)+c+t)}
\end{equation}
for $1\le t\le\min{\{b,p-c\}}$ and
\begin{equation}\label{eq:A2}
\mathscr{P}^{p^at}QA^{2p^a(pb+c)}=0
\qquad
\text{in $QA^{2p^a(p(b+t)+c-t)}$}
\end{equation}
for $c\le t<p$.

\textup{(2)}\qua If $a\ge 0$ and $0<c<p$,
then
\begin{equation}\label{eq:A3}
\mathscr{P}^{p^at}\co QA^{2p^ac}
\longrightarrow QA^{2p^a(tp+c-t)}
\end{equation}
is an isomorphism for $1\le t<c$.
\end{prop}

\begin{proof}
First we consider the case of $a=0$.
Let us prove (1) by downward induction on $b$.
If $b$ is large enough,
then the result is clear since $QA^{2(pb+c)}=0$.
Assume that $y_j$ is one of the generators
of \eqref{eq:truncated} for $1\le j\le l$
and $\deg y_j=2(pb+c)$ with $b>0$ and $0<c<p$.
By inductive hypothesis,
we can assume that if $f>b$ and $0<g<p$,
then
\begin{equation}\label{eq:induction}
QA^{2(pf+g)}=\mathscr{P}^tQA^{2(p(f-1)+g+1)}
\end{equation}
for $1\le t\le\min{\{f,p-g\}}$.
If we put
\begin{equation*}
\beta=\frac{1}{pb+c}\mathscr{P}^{pb+c-1}(y_j)\in A^{2(p(pb+c-1)+1)},
\end{equation*}
then by \eqref{eq:induction},
we have
\begin{equation*}
\beta-\mathscr{P}^{p-1}(\gamma)\in DA^*
\end{equation*}
for some $\gamma\in QA^{2p(p(b-1)+1)}$.
Since $A^*$ is a $\mathscr{D}_n$--algebra,
\begin{equation}\label{eq:beta}
\beta-\frac{1}{pb+c}\mathscr{P}^{pb+c-1}(\mu)
-\mathscr{P}^{p-1}(\gamma-\nu)\in D^{n+1}A^*
\end{equation}
for some decomposable classes
$\mu\in DA^{2(pb+c)}$
and $\nu\in DA^{2p(p(b-1)+1)}$.
If we apply $\mathscr{P}^1$ to \eqref{eq:beta},
then $y_j^p=\mathscr{P}^1(\xi)$
for some $\xi\in D^{n+1}A^*$
since $\mathscr{P}^{pb+c}(\mu)=\mu^p=0$
in $A^*$ and $\mathscr{P}^1\mathscr{P}^{p-1}=p\mathscr{P}^p=0$.
Then for some generator $y_i$, $\mathscr{P}^1(y_i)$ must contain
some $y_j^t$ with $1\le t\le p$ and $t+n=p$.
By the assumption of $n>(p-1)/2$,
we have $1\le t\le n$,
which implies that $y_j=\mathscr{P}^1(y_k)$
for some $1\le k\le l$
by \fullref{prop:propthmAB}.
By iterating this argument,
we have \eqref{eq:A1}.

Now \eqref{eq:A2}
follows from \eqref{eq:A1}.
In fact, if $y_j$ is a generator in \eqref{eq:truncated} with
$\deg y_j=2(pb+c)$ for some $b>0$ and $0<c<p$, then
we show that $\mathscr{P}^c(y_j)=0$.
If $b+c<p$,
then by \eqref{eq:A1},
we have $y_j=\mathscr{P}^b(\kappa)$
for $\kappa\in QA^{2(b+c)}$,
which implies that
\begin{equation}\label{eq:kappa}
\mathscr{P}^c(y_j)=\mathscr{P}^c\mathscr{P}^b(\kappa)
={\binom{b+c}{b}}\kappa^p=0
\end{equation}
in $QA^{2p(b+c)}$.
On the other hand,
if $p\le b+c$,
then by \eqref{eq:A1},
we have $y_j=\mathscr{P}^{p-c}(\zeta)$
for $\zeta\in QA^{2p(b+c-p+1)}$,
and so
\begin{equation}\label{eq:zeta}
\mathscr{P}^c(y_j)=\mathscr{P}^c\mathscr{P}^{p-c}(\zeta)
={\binom{p}{c}}\mathscr{P}^p(\zeta)=0.
\end{equation}

Next we show (2) with $a=0$.
We only have to show that $\mathscr P^{c-1}$ is a monomorphism
on $QA^{2c}$.
Let $y_j$ be a generator
in \eqref{eq:truncated} such that $\deg y_j=2c$ with $0<c<p$.
Suppose contrarily that $\mathscr{P}^{c-1}(y_j)=0$ in $QA^{2(c-1)p+1}$.
Since $A^*$ is a $\mathscr{D}_n$--algebra,
we have that
\begin{equation*}
\mathscr{P}^{c-1}(y_j-\mu)\in D^{n+1}A^{2(c-1)p+1}
\end{equation*}
for some decomposable class
$\mu\in DA^{2c}$.
Then by a similar argument to the proof of (1),
we have that $y_j=\mathscr{P}^1(y_k)$
for some $1\le k\le l$ with $\deg y_k=2(c-p+1)$,
which is impossible for dimensional reasons.
This completes the proof
of \fullref{prop:propthmA}
in the case of $a=0$.

Let $I$ denote the ideal of $A^*$
generated by $y_i$
with $m_i\not\equiv 0$ mod $p$.
Then for dimensional reasons and
by \eqref{eq:kappa} and \eqref{eq:zeta},
we see that $I$ is closed
under the action of $\mathscr{A}_p^*$,
which implies that $A^*/I$
is an unstable $\mathscr{A}_p^*$--algebra
given by
\begin{equation*}
A^*/I=T^{[p+1]}[y_{i_1},\ldots,y_{i_q}]
\qquad
\text{with $m_{i_d}\equiv 0 \bmod p$}
\end{equation*}
for $1\le d\le q$.
Set $m_{i_d}=ph_d$ with $h_d\ge 1$
for $1\le d\le q$.
Let $B^*$ denote the truncated polynomial algebra
at height $p+1$ given by
\begin{equation*}
B^*=T^{[p+1]}[z_1,\ldots,z_q]
\qquad
\text{with $\deg z_d=h_d$}
\end{equation*}
for $1\le d\le q$.
If we define a map
$\widetilde{\mathscr{L}}\co\{y_{i_1},\ldots,y_{i_q}\}\to B^*$
by $\widetilde{\mathscr{L}}(y_{i_d})=z_d$ for $1\le d\le q$,
then $\widetilde{\mathscr{L}}$ is extended to
an isomorphism $\mathscr{L}\co A^*/I\to B^*$.
Moreover, $B^*$ admits an unstable
$\mathscr{A}_p^*$--algebra structure
by the action $\mathscr{P}^r(z_d)
=\mathscr{L}(\mathscr{P}^{pr}(y_{i_d}))$
for $r\ge 1$.
Then we can show that $B^*$ is a $\mathscr{D}_n$--algebra
concerning this structure since so is $A^*$.
From the above arguments,
we have the required results for $B^*$
in the case of $a=0$,
which implies that $A^*$
satisfies the required results for $a=1$.
By repeating these arguments,
we can show that $A^*$
satisfies the desired conclusions
of \fullref{prop:propthmA}
for any $a\ge 0$.
This completes the proof.
\end{proof}

Now we prove \fullref{thmA} as follows:

\begin{proof}[Proof of~\fullref{thmA}]
By Browder~\cite[Theorem~8.6]{browder},
$H^*(X;\mathbb{Z}/p)$,
the mod $p$ cohomology,
is an exterior algebra in \eqref{eq:exterior}.
Let $\widetilde{X}$ be the universal cover of $X$.
From the proof of \cite[Lemma~3.9]{hemmi-kawamoto},
we have that $\widetilde{X}$
is a simply connected $A_p$--space
admitting an $AC_n$--form.
It is enough to prove \fullref{thmA} for $\widetilde{X}$
since $X\simeq\widetilde{X}\times T$
for a torus $T$ by Kane~\cite[page~24]{kane}.
By \fullref{thm:acn-quasicn}
and \fullref{thm:dn-algebra},
we have that $A^*(\widetilde{X})$
is a $\mathscr{D}_n$--algebra.
Then by \fullref{thm:modify} (3)
and \fullref{prop:propthmA},
we have the required conclusion.
This completes the proof of \fullref{thmA}.
\end{proof}

By using \fullref{thmA},
we prove \fullref{thmB} as follows:

\begin{proof}[Proof of \fullref{thmB}]
We proceed by using a similar way
to the proof of \cite[Theorem~1.1]{hemmi2}.
Since $X\simeq\widetilde{X}\times T$
for a torus $T$ as in the proof of \fullref{thmA},
the Steenrod operations $\mathscr{P}^j$
act trivially on $QH^*(\widetilde{X};\mathbb{Z}/p)$
for $j\ge 1$.
By \fullref{thmA},
if $QH^{2m-1}(\widetilde{X};\mathbb{Z}/p)\ne 0$,
then $m=p^a$ for some $a\ge 1$,
and so the mod $p$ cohomology of $\widetilde{X}$
is an exterior algebra in \eqref{eq:exterior},
where $m_i=p^{a_i}$ with $a_i\ge 1$
for $1\le i\le l$.

Let $P_{p-1}(\widetilde{X})$ be
the $(p-1)$-th projective space of $\widetilde{X}$.
Then by \eqref{eq:cohomology of projective},
there is an ideal $S_{p-1}
\subset H^*(P_{p-1}(\widetilde{X});\mathbb{Z}/p)$
closed under the action of $\mathscr{A}_p^*$
with
\begin{equation*}
H^*(P_{p-1}(\widetilde{X});\mathbb{Z}/p)/S_{p-1}
\cong T^{[p]}[y_1,\ldots,y_l],
\end{equation*}
where $T^{[p]}[y_1,\ldots,y_l]$ is
the truncated polynomial algebra at height $p$
generated by $y_i\in H^{2p^{a_i}}(P_{p-1}(\widetilde{X});\mathbb{Z}/p)$
with $\iota_1^*\ldots\iota_{p-2}^*(y_i)=\sigma(x_i)
\in H^{2p^{a_i}}(\Sigma\widetilde{X};\mathbb{Z}/p)$.
Moreover, we have that the composite
\begin{equation}\label{eq:iso-epi}
\begin{CD}
H^t(P_p(\widetilde{X});\mathbb{Z}/p)
@>{\iota_{p-1}^*}>> H^t(P_{p-1}(\widetilde{X});\mathbb{Z}/p)
@>>> T^{[p]}[y_1,\ldots,y_l]
\end{CD}
\end{equation}
is an isomorphism for $t<2p^{a_1+1}$
and an epimorphism for $t<2(p^{a_1+1}+p^{a_1}-1)$
by \cite[page~106, (4.10)]{hemmi2}.
As in \cite[page~106, (4.11)]{hemmi2},
we can show
\begin{equation}\label{eq:im-p-intersection}
\textup{Im}\ \mathscr{P}^{p^{a_1}}\cap H^t(P_p(\widetilde{X});\mathbb{Z}/p)=0
\end{equation}
for $t\le 2p^{a_1+1}$.
In fact,
by \eqref{eq:iso-epi} and for dimensional reasons,
we have
\begin{equation*}
\begin{split}
\textup{Im}\ \beta\cap H^t(P_p(\widetilde{X});\mathbb{Z}/p) & = 0 \\
\textup{Im}\ \mathscr{P}^1\cap H^t(P_p(\widetilde{X});\mathbb{Z}/p) & = 0
\end{split}
\end{equation*}
for $t\le 2p^{a_1+1}$,
which implies \eqref{eq:im-p-intersection}
by Liulevicius~\cite{liulevicius}
or Shimada--Yamanoshita~\cite{shimada-yamanoshita}.
Now we can choose $w_1 \in H^{2p^{a_1}}(P_p(\widetilde{X});\mathbb{Z}/p)$
with $\smash{\iota_{p-1}^*(w_1)}=y_1$ by \eqref{eq:iso-epi},
and we have $\smash{w_1^p}=\mathscr{P}^{p^{a_1}}(w_1)=0$
by \eqref{eq:im-p-intersection}.
Then we have a contradiction
by using the same argument
as the proof of \cite[Theorem~1.1]{hemmi2},
and so $\widetilde{X}$ is contractible,
which implies that $X$ is a torus.
This completes the proof of \fullref{thmB}.
\end{proof}

To show \fullref{thmC},
we need the following definition:

\begin{definition}\label{def:acnform-in-x}
Assume that $X$ is an $A_n$--space
and $Y$ is a space.

\textup{(1)}\qua
An $AC_n$--form on a map $\phi\co Y\to X$
is a family of maps
$\{R_i\co\Gamma_i\times Y^i\to X\}_{1\le i\le n}$
with the conditions $R_1(*,y)=\phi(y)$ for $y\in Y$
and \eqref{eqn:acnsp2}--\eqref{eqn:acnsp3}.

\textup{(2)}\qua
A quasi $C_n$--form on a map
$\kappa\co\Sigma Y\to\Sigma X$
is a family of maps
$\{\zeta_i\co J_i(\Sigma Y)\to P_i(X)\}_{1\le i\le n}$
with the conditions $\zeta_1=\kappa$
and \eqref{eq:quasicn2}.
\end{definition}

By using the same argument as the proof
of \fullref{thm:acn-quasicn} (1),
we can prove the following result:

\begin{theorem}\label{thm:acn-quasicn-map}
Assume that X is an $A_n$--space,
$Y$ is a space and $\phi\co Y\to X$ is a map.
Then any $AC_n$--form on $\phi$ induces
a quasi $C_n$--form on $\Sigma\phi$.
\end{theorem}

Now we prove \fullref{thmC} as follows:

\begin{proof}[Proof of \fullref{thmC}]
First we show that if $X$ admits an $AC_n$--form,
then $nm_l\le p$.

We prove by induction on $n$.
If $n=1$,
then the result is proved
by Hubbuck--Mimura~\cite{hubbuck-mimura}
and Iwase~\cite[Proposition~0.7]{iwase2}.
Assume that
the result is true for $n-1$.
Then by inductive hypothesis,
we have $(n-1)m_l\le p$.
Now we assume that $X$
admits an $AC_n$--form with
\begin{equation}\label{eq:induction2}
(n-1)m_l\le p<nm_l.
\end{equation}
Then we show a contradiction.

Let $\widetilde{X}$ be the universal
covering space of $X$.
Then $\widetilde{X}$ is a simply connected
$A_p$--space mod $p$ homotopy equivalent to
\begin{equation}\label{eq:product of spheres}
S^{2m_1-1}\times \cdots \times S^{2m_l-1}
\qquad
\text{with $1< m_1\le \cdots \le m_l$}
\end{equation}
and the multiplication of $\widetilde{X}$
admits an $AC_n$--form
by \cite[Lemma~3.9]{hemmi-kawamoto}.
Now we can set that
\begin{equation*}
A^*(\widetilde{X})=T^{[p+1]}[y_1,\ldots,y_l]
\qquad
\text{with $\deg y_i=2m_i$}
\end{equation*}
for $1\le i\le l$,
where $1<m_1\le\cdots\le m_l\le p$.
By \fullref{thm:acn-quasicn}
and \fullref{thm:dn-algebra},
$A^*(\widetilde{X})$ is a $\mathscr{D}_n$--algebra.

First we consider the case of $m_l<p$.
Let $J$ be the ideal of $A^*(\widetilde{X})$
generated by $y_i$ for $1\le i\le l-1$.
Then we see that
\begin{equation}\label{eq:idealj}
\mathscr{P}^1(y_i)\not\in J
\qquad
\text{for some $1\le i\le l$.}
\end{equation}
In fact,
if we assume that $\mathscr{P}^1(y_i)\in J$
for any $1\le i\le l$,
then $\mathscr{P}^1(y_l)\in J$ and $\mathscr{P}^1(J)\subset J$.
This implies that
\begin{equation*}
y_l^p=\mathscr{P}^{m_l}(y_l)=\frac{1}{m_l!}(\mathscr{P}^1)^{m_l}(y_l)\in J,
\end{equation*}
which is a contradiction,
and so we have \eqref{eq:idealj}.
Then for dimensional reasons
and by \eqref{eq:induction2},
\begin{equation*}
2(n-1)m_l<\deg\mathscr{P}^1(y_i)<2(n+1)m_l,
\end{equation*}
which implies that $\mathscr{P}^1(y_i)$ contains the term
$ay_l^n$ with $a\ne 0$ in $\mathbb{Z}/p$
by \eqref{eq:idealj}.
By \fullref{prop:propthmAB},
we have $y_l\in\mathscr{P}^1QA^{2(m_l-p+1)}(\widetilde{X})$,
which causes a contradiction since $m_l<p$.

Next let us consider the case of $m_l=p$.
In this case, \eqref{eq:induction2} is equivalent to $n=2$,
and so $\widetilde{X}$ is assumed to have an $AC_2$--form.
Then from the same arguments as above,
we have that $A^*(\widetilde{X})$
is a $\mathscr{D}_2$--algebra.
By Kanemoto~\cite[Lemma~3]{kanemoto},
there is a generator $y_k\in QA^{2(p-1)}(\widetilde{X})$
for some $1\le k<l$.
Let $K$ be the ideal of $A^*(\widetilde{X})$
generated by $y_i$ with $i\ne k$.
From the same reason as \eqref{eq:idealj},
we see that $\mathscr{P}^1(y_i)\not\in K$
for some $1\le i\le l$.
Then for dimensional reasons,
we see that $\mathscr{P}^1(y_i)$ contains
the term $by_k^2$ with $b\ne 0$ in $\mathbb{Z}/p$.
By \fullref{prop:propthmAB},
we have a contradiction,
and so $\widetilde{X}$ does not admit an $AC_2$--form.

Next we show that if $nm_l\le p$,
then $X$ admits an $AC_n$--form.
Since it is clear for $n=1$ or $m_l=1$,
we can assume that $nm_l<p$.
Let $Y$ denote the wedge sum of spheres
given by
\begin{equation*}
Y=(S^{2m_1-1}\vee\ldots\vee S^{2m_l-1})_{(p)}
\end{equation*}
with the inclusion $\phi\co Y\to X$.
First we construct an $AC_n$--form
$\{R_i\co\Gamma_i\times Y^i\to X\}_{1\le i\le n}$
on $\phi\co Y\to X$.

Suppose inductively that
$\{R_i\}_{1\le i< t}$ are constructed for some $t\le n$.
Then the obstructions for the existence
of $R_t$ belong to the following
cohomology groups for $j\ge 1$:
\begin{equation}\label{eq:obstruction}
H^{j+1}(\Gamma_t\times Y^t,
\partial\Gamma_t\times Y^t\cup\Gamma_t\times Y^{[t]};
\pi_j(X))
\cong\widetilde{H}^{j+2}((\Sigma Y)^{(t)};\pi_j(X))
\end{equation}
since $\Gamma_t\times Y^t/(\partial\Gamma_t\times Y^t
\cup\Gamma_t\times Y^{[t]})
\simeq\Sigma^{t-1}Y^{(t)}$.
This implies that \eqref{eq:obstruction}
is non-trivial only if 
$j$ is an even integer
with $j<2p-2$ since
\begin{equation*}
\Sigma Y\simeq(S^{2m_1}\vee\ldots\vee S^{2m_l})_{(p)}
\end{equation*}
and $tm_l \le nm_l<p$.
On the other hand,
according to Toda~\cite[Theorem~13.4]{toda},
$\pi_j(X)= 0$ for any even integer $j$ with $j< 2p-2$
since $X$ is given by \eqref{eq:product of spheres}.
Thus \eqref{eq:obstruction} is trivial for all $j$, and we
have a map $R_t$.
This completes the induction, and 
we have an $AC_n$--form $\{R_i\}_{1\le i\le n}$
on $\phi\co Y\to X$.

Since $X$ is an $H$--space,
there is a map $\beta\co\Omega\Sigma X\to X$
with $\beta\alpha\simeq 1_{X}$,
where $\alpha\co X\to\Omega\Sigma X$
denotes the adjoint of $1_{\Sigma X}\co\Sigma X\to\Sigma X$.
Moreover,
$\Sigma Y$ is a retract of $\Sigma X$,
and so we have a map $\nu\co\Sigma X\to\Sigma Y$
with $\nu(\Sigma\phi)\simeq 1_{\Sigma Y}$.
Put $\lambda=\beta\theta\co X\to X$,
where $\theta\co X\to\Omega\Sigma X$
denotes the adjoint of $(\Sigma\phi)\nu$.
Then we see that $\lambda$ induces
an isomorphism on the mod $p$ cohomology,
and so $\lambda$ is a mod $p$ homotopy equivalence.

By \fullref{thm:acn-quasicn-map},
there is a quasi $C_n$--form
$\{\zeta_i\co J_i(\Sigma Y)\to P_i(X)\}_{1\le i\le n}$
on $\Sigma\phi\co\Sigma Y\to\Sigma X$.
Let $\xi_i\co J_i(\Sigma X)\to P_i(X)$
be the map defined by
$\xi_i=\zeta_iJ_i(\nu(\Sigma\lambda^{-1}))$
for $1\le i\le n$,
where $\lambda^{-1}\co X\to X$ denotes
the homotopy inverse of $\lambda$.
Then the family $\{\xi_i\}_{1\le i\le n}$
satisfies that $\xi_i|_{J_{i-1}(\Sigma X)}
=\iota_{i-1}\xi_{i-1}$ for $2\le i\le n$
and $\xi_1=(\Sigma\phi)\nu(\Sigma\lambda^{-1})
=\chi(\Sigma\theta)(\Sigma\lambda^{-1})$,
where $\chi\co\Sigma\Omega\Sigma X\to \Sigma X$
is the evaluation map.
Since $\iota_1(\Sigma\beta)\simeq\iota_1\chi
\co\Sigma\Omega\Sigma X\to P_2(X)$
by Hemmi~\cite[Lemma~2.1]{hemmi4},
we have $\xi_2|_{\Sigma X}=\iota_1\xi_1
\simeq\iota_1$.
Let $\psi_i\co J_i(\Sigma X)\to P_i(X)$
be the map defined by
$\psi_1=1_{\Sigma X}$ and $\psi_i=\xi_i$
for $2\le i\le n$.
Then the family $\{\psi_i\}_{1\le i\le n}$
satisfies \eqref{eq:quasicn1}--\eqref{eq:quasicn2}.
By \fullref{thm:acn-quasicn} (2)
and \fullref{rem:acn-quasicn},
we have an $AC_n$--form
$\{Q_i\co\Gamma_i\times X^i\to X\}_{1\le i\le n}$
on $X$ with \eqref{eqn:acnsp1}--\eqref{eqn:acnsp3}.
This completes the proof of \fullref{thmC}.
\end{proof}

\bibliographystyle{gtart}
\bibliography{link}

\end{document}